\documentclass[12pt]{article}

\newcommand{\cv}{{\cal V}}
\usepackage{amsmath}

\input{diagrams}
\newenvironment{ting}[1]%
{%
\begin{list}{\textbf{#1} }%
{\setlength{\labelsep}{0mm}\setlength{\leftmargin}{0mm}%
\setlength{\labelwidth}{0mm}\setlength{\listparindent}{\parindent}%
\setlength{\parsep}{\parskip}\setlength{\partopsep}{0mm}}%
\item \it
}%
{%
\end{list}%
}
\begin{document}
\title{Commutation Structures}
\author{Anders Kock}
\date{}
\maketitle

\small  \quad \begin{minipage}{120mm}{\bf Abstract} For fixed object $X$ in a monoidal category, an 
$X$-commutation structure on an object $A$ is just a map 
$X\otimes A \to A \otimes X$. We study aspects of such structures in 
case $A$ has a dual object.\end{minipage}

\bigskip
\bigskip
\normalsize

We consider a monoidal category $\cv , \otimes, I$; for simplicity 
we let it be strict (the application we have in mind is anyway a 
category of endofuncors on a category, with composition as $\otimes$).

Let $X$ be an object in $\cv$, fixed throughout. An {\em 
$X$-commutation structure} on an object $A\in \cv$ is an arrow
$$\alpha : X\otimes A \to A\otimes X.$$ A {\em morphism} of 
$X$-commutation structures $(A,\alpha ) \to (B, \beta )$ is an arrow 
$A\to B$ such that the obvious square
$$\begin{diagram}
X\otimes A &\rTo ^{\alpha} & A\otimes X\\
\dTo ^{X\otimes f}&& \dTo _{f\otimes X}\\
X \otimes B &\rTo _{\beta }& B\otimes X
\end{diagram}$$
commutes.

In this way, we get a category of $X$-commutation structures; denote 
it $X\mbox{-}\cv $. There is a faithful forgetful functor $X\mbox{-}\cv  \to \cv$:
to $(A, \alpha )$, associate $A$.

There is also a monoidal structure on $X\mbox{-}\cv $, preserved strictly by the 
forgetful functor $X\mbox{-}\cv  \to \cv$: 
If $(A,\alpha )$ and $(B, \beta )$ are objects in 
$X\mbox{-}\cv $, we get a commutation structure $\gamma$ on $A\otimes B$ in an obvious 
way: $\gamma$ is taken to be the composite
$$\begin{diagram}
X \otimes A \otimes B & \rTo ^{\alpha \otimes B} & A \otimes X \otimes 
B & \rTo ^{A \otimes \beta } & A \otimes B \otimes X .
\end{diagram}$$
Also $I$ carries a canonical $X$-commutation structure, namely the 
canonical $X\otimes I \to I \otimes X$ (an identity, in fact, since 
we assumed $\cv$ strict).

Recall that a right dual for an object $A$ in a monoidal category is 
an object $B$ together with arrows (``unit and counit'')
$$ \eta: I \to 
B\otimes A\quad \mbox{ and } \quad \epsilon : A\otimes B \to I$$
satisfying the usual two triangle equations.

We are interested in the monoidal category $X\mbox{-}\cv $:

\begin{ting}{Theorem 1}If an object $(A, \alpha ) \in X\mbox{-}\cv $ admits a right dual, 
then $\alpha$ is an invertible arrow in $\cv$.
\end{ting}
{\bf Proof.} The assumption on $(A, \alpha )$ can be expressed: $A$ 
admits a right dual $B$ in the category $\cv$, and there is an 
$X$-commutation structure $\beta$ on $B$ in such a way that the unit 
and counits are morphisms of commutation structures. 

The proof is now purely equational: we exhibit a two sided inverse 
for $\alpha : X\otimes A \to A \otimes X$.
We shall prove that the following composite $\gamma$ will serve. (To save 
space, we write $\otimes$ just as a dot, $A\cdot B$ for $A \otimes 
B$, etc.)
$$\begin{diagram}
A\cdot X & \rTo ^{A\cdot X \cdot \eta}&A \cdot X \cdot B 
\cdot A & \rTo ^{A\cdot \beta \cdot A} & A\cdot B \cdot X 
\cdot A & \rTo ^{\epsilon \cdot X \cdot A} & X \cdot A .
\end{diagram}$$
To prove that $\alpha \circ \gamma$ is the identity means to prove 
that the clockwise composite in the following diagram is the 
identity, and this is proved by considering the rest of the diagram, 
as we shall argue:
$$\begin{diagram}[width=5em, nohug]
A\cdot X & \rTo ^{A\cdot X \cdot \eta}&A \cdot X \cdot B 
\cdot A & \rTo ^{A\cdot \beta \cdot A} & A\cdot B \cdot X 
\cdot A & \rTo ^{\epsilon \cdot X \cdot A} & X \cdot A \\
&\rdTo(4,2) _{A\cdot \eta \cdot X }&&&\dTo _{A\cdot B \cdot \alpha}&&\dTo 
_{\alpha}\\
&&&&A\cdot B \cdot A \cdot X& \rTo _{\epsilon \cdot A \cdot X} & A 
\cdot X
\end{diagram}$$
Here the square commutes, because $\otimes$ is a functor in two 
variables. For the left hand cell, consider the square (with top map 
the identity)
\begin{equation}\begin{diagram}
X\cdot I & \rTo   &&& I\cdot X\\
\dTo ^{X\cdot \eta}&&&&
\dTo _{\eta \cdot X}\\
X\cdot B\cdot A & \rTo _{\beta \cdot A}&B\cdot X \cdot A& \rTo _{B \cdot\alpha }& 
B\cdot A \cdot X.
\end{diagram}\label{one}\end{equation}
It commutes by the assumption that $\eta : I \to B\cdot A$ is a 
morphism of $X$-commutation structures. If we apply the functor $A 
\cdot -$, we get (by a geometric reflection) the desired 
commutativity of left hand cell. -- Finally, the lower 
(counter-clockwise) composite is an identity arrow: it is just the 
functor $_\cdot X$ applied to one of the triangle equations for 
$\eta ,\epsilon$. Thus $\alpha \circ \gamma$ is the identity arrow of 
$A\cdot X$.

The proof that $\gamma \circ \alpha$ is the identity is much similar. 
The map $\gamma \circ \alpha $ appears as the counterclockwise composite in 
the diagram
$$\begin{diagram}[width=5em]
X\cdot A & \rTo ^{X\cdot A \cdot \eta}& X\cdot A \cdot B \cdot A &&&&\\
\dTo ^{\alpha}&& \dTo _{\alpha \cdot B \cdot A}&\rdTo(4,2) ^{X\cdot 
\epsilon \cdot A}&&&\\
A\cdot X & \rTo _{A\cdot X \cdot \eta }& A\cdot X \cdot B \cdot A & 
\rTo _{A\cdot \beta \cdot A }& A\cdot B \cdot X \cdot A & \rTo 
_{\epsilon \cdot X \cdot A}& X\cdot A.
\end{diagram}$$
The square on the left commutes because $\otimes$ is a functor in two 
variables, and the cell on the right commutes because $\epsilon$ was 
assumed to be a morphism of $X$-commutation structures (we omit the 
diagram, which is analogous to (\ref{one})). Finally, the upper 
(clockwise) composite is the identity: it is just the functor $- \cdot 
A$ applied to one of the triangle equations for $\eta, \epsilon$. 
Thus $\gamma \circ \alpha$ is the identity arrow of $X \cdot A$. This 
proves the Theorem.

\medskip
{\bf Example.} Let ${\cal C}$ be a category with coproducts. Let 
$\cv$ be the monoidal category of endofunctors on ${\cal C}$. Let $J$ 
be  a fixed set, and let $X$ be the endofunctor on ${\cal C}$ given 
by $C \mapsto \coprod_J C$ (coproduct of $J$ copies of $C$). Let $A: 
{\cal C}\to {\cal C}$ be any endofunctor. There is 
for each $C \in {\cal C}$ a map
$$\coprod _J A(C) \to A(\coprod _J C)$$ which on the $j$'th summand 
of $\coprod _J A(C)$ returns $A(\mbox{incl}_j ): A(C) \to A(\coprod 
_J C )$. This is natural in $C$, and thus is a natural transformation 
$$\alpha : X \circ A \to A \circ X ,$$
in other words, an $X$-commutation structure on $A$ in $\cv$. It is 
clear that any natural transformation $A_1 \to A_2$ between 
endofunctors on ${\cal C}$ is a morphism of $X$-commutations. In 
particular, if $A$ has a right adjoint ( = a right dual in the 
monoidal category $\cv$ of endofunctors), this adjointness (duality) lifts to 
an adjointness/duality in $X\mbox{-}\cv $.

From the Theorem then follows that $\alpha$ is actually an isomorphism. 
This can be restated: ``If $A$ has a right adjoint, then 
$A$ commutes with copowers''; which is of course no big surprise.

More generally, if ${\cal C}$ is a category enriched over a category 
${\cal S}$, it makes sense to say that it is {\em tensored} over ${\cal 
S}$. To say that an endofunctor $A$ on ${\cal C}$ is {\em enriched} 
can expressed in terms of existence of {\em a tensorial strength} 
$\alpha$, cf.\ \cite{SFMM}, which is a map, natural in ($J$ and) $C$,
$$J\otimes A(C) \to A(J \otimes C)$$
for $J \in {\cal S}$ and $C\in {\cal C}$. Now an adjointness $A \dashv 
B$ is no longer automatically enriched/strong, but if it is, the 
Theorem implies that $A$ commutes with tensors $J\otimes -$ up to 
isomorphism.  

\medskip

Consider a Cartesian Closed Category ${\cal S}$. Being closed, it is 
enriched over itself, and the tensors $J\otimes C$ are just $J \times 
C$. An enrichment/strength of an endofunctor $A: {\cal S} \to {\cal 
S}$ then can be encoded as a ``tensorial strength'', i.e.\
as a natural family of maps
$$J\times A(C) \to A(J\times C ),$$
equivalently as a commutation (for each $J$)
$$\alpha :X\circ A \to A\circ X$$
where $X$ denotes the endofunctor $J\times - $.

Consider in particular an  endofunctor of the form 
$(-)^D : {\cal S} \to {\cal S}$. If this endofunctor has a right 
adjoint, $D$ is called an {\em atom}, cf.\ e.g. \cite{SDG}, 
\cite{AFEF}, or \cite{Yetter} (who calls such objects $D$ {\em tiny}).
A Corollary of the Theorem is now a result (due to Yetter):

\begin{ting}{Propoition 2}If $(-)^D$ has a {\em strong} right adjoint, then $D=1$.
\end{ting}
{\bf Proof.} Let $X$ be an endofunctor of the form $J \times -$ ($J$ 
an arbitrary object of ${\cal S}$). If $(-)^D$ has a strong right 
adjoint, the Theorem implies that the 
natural map (the tensorial strength, or commutation structure)
$$J\times Y^D \to (J\times Y)^D$$
is an isomorphism for any $Y$. This in particular applies to $Y=1$, so that
$J \times 1^D \cong (J \times 1)^D$, but since $1^D \cong 1$, we get 
from this that the natural map
$$J \to J^D$$
is an isomorphism, for any $J$. A suitable enriched Yoneda Lemma now 
gives the result, but there is an elementary proof: 

From the fact that the natural map $J \to J^D$ is an isomorphism, we 
conclude that for each $X$, there is  a bijection 
between $\hom (X,J)$ and $\hom (X,J^D )$ (induced by the natural map
$J \to J^D$); passing to transposes, there 
is a bijection between $\hom (X, J)$ and $\hom(X\times D ,J)$ (induced 
by the projection $X\times D \to X $). Now take $X=1$ and $J=D$ and 
conclude that the projection $1\times D \to D$ factors across the 
projection $1\times D \to 1$. From this follows that $D$ is a retract 
of $1$, hence is itself (isomorphic to ) $1$.

 %
\medskip

The notion of $X$-co-commutation is obtained by duality, thus a 
$X$-co-commutation on $A$ is a map $a:A\cdot X \to X \cdot A$. 
Similarly for morphisms of co-commutations $(A,a)\to (B,b)$,
and for the monoidal category 
of co-commutations. In particular, if $A \dashv B$ by virtue of 
$\eta , \epsilon$, as above, and these are compatible with the 
co-commutations, it follows by Theorem 1 (dualized) that the 
co-commutation $b$ on $B$ is invertible. Thus, the $b^{-1}$ in the 
following Theorem makes sense. (I omit the $\otimes$, formerly 
abbreviated $\cdot$; now 
they are both denoted just by concatenation.)

\begin{ting}{Theorem 3} Let $(A,a) \dashv (B, b)$ in the category of 
$X$-co-commutations. Then $a$ can be expressed in terms of $b^{-1}$ as 
follows: $a=$
\begin{equation}\begin{diagram} A X & \rTo^{AX\eta} & AXBA &\rTo^{Ab^{-1}A}& ABXA &
\rTo^{\epsilon 
XA} & XA .\end{diagram}\label{two}\end{equation}
\end{ting}
{\bf Proof.} Let us denote by $c$ the composite co-commutation $BAX 
\to XBA$ on $BA$, Thus the middle of the small triangles in the 
following diagram commutes by definition:
$$\begin{diagram}[nohug]
AX & \rTo ^{AX\eta} & AXBA &&\\
\dTo ^{A\eta X}& \ruTo ^{Ac} & \uTo _{AbA}&\rdTo ^{Ab^{-1}A}&\\
ABAX & \rTo _{ABa}& ABXA & \rTo _{=} & ABXA\\
\dTo ^{\epsilon AX}&&&&\dTo _{\epsilon XA}\\
AX &&\rTo _a && XA
\end{diagram}$$
The upper left little triangle commutes since $\eta : I \to BA$ is a 
morphism of co-commutation structures. The next triangle commutes by 
definition of $c$, as observed, and the third triangle evidently 
commutes. The bottom ``square'' commutes by bi-functorality of 
$\otimes$. finally, the left hand column is the identity map of $AX$, 
by one of the triangle equations for $\eta$, $\epsilon$. Thus, the 
counterclockwise composite in the diagram is $a$, the clockwise is the 
arrow in (\ref{two}). This proves the Theorem.

\medskip

My motivation for this Theorem was the desire to understand the proof 
of Proposition XIV.3.1 in \cite{Kassel}; according to this, we have in 
a braided category that the braiding $c_{AX} : AX \to XA$
of a dualizable object $A$ is determined by  the 
braiding $c_{BX}$ for its dual $B$ ($X$ a fixed object). (For the 
comparison with Kassel, our $A$ is his 
$V^*$, $B$ is $V$, and $X$ is $W$ .)

More precisely, let $X$ be a fixed object in a braided monoidal 
category. Then for any object $D$, the braiding
$c_{DX}$ defines an $X$- co-commutation structure on $D$, and any map 
 $D_1 \to D_2$ is a homomorphism of 
co-commutativity structures, just by naturality of $c_{-,X}$. Also, 
the composite co-commutativity structure on $BA$ is $c_{BA,X}$, by 
the ``hexagon'' axiom for braidings (which here reduces to a 
triangle), and any duality $A\dashv B$ is automatically compatible 
with the co-commutations (since any map is).

Applying Theorem 3 to $a := c_{A,X}$, $b:= c_{B,X}$ therefore 
expresses $c_{A,X}$ in terms of $c_{B,X}$, as in the statement of the 
Proposition in \cite{Kassel} (the proof in loc.\ cit.\ is given in 
terms of ``graphical calculus'').

\end{document}